\documentclass[english]{article}
\usepackage[T1]{fontenc}
\usepackage[latin9]{inputenc}
\usepackage{amsmath}
\usepackage{amsthm}
\usepackage{verbatim}
\usepackage{amssymb}

\makeatletter
\theoremstyle{plain}
\ifx\thechapter\undefined
\newtheorem{thm}{\protect\theoremname}
\else
\newtheorem{thm}{\protect\theoremname}[chapter]
\fi
\theoremstyle{plain}
\newtheorem*{lem*}{\protect\lemmaname}
\theoremstyle{remark}
\newtheorem*{rem*}{\protect\remarkname}

\usepackage{geometry}

\makeatother

\usepackage{babel}
\providecommand{\lemmaname}{Lemma}
\providecommand{\remarkname}{Remark}
\providecommand{\theoremname}{Theorem}

\begin{document}
	\title{ On the problem of Kakeya}
	\author{Rado\v s Baki\'{c}}
	\date{}
	\maketitle
	\begin{abstract} 
		Let $p(z)$ be a complex polynomial of degree $n$, having $k$ its zeros in the unite disc. We prove that at least one zero of $p^{(k-1)}(z)$ lie in the disc $|z|\le \frac{2(n-k+1)}{\ln 2} $.

\end{abstract}
	\textbf{Key words:}  zeros of polynomial, critical points of  a polynomial, apolar polynomials.
	\\
	\textbf{AMS Subject Classification} Primary 26C10, Secondary 30C15.

	\vspace{1cm}
	
 	Let $p(z)$ be a complex polynomial of degree $n$, and suppose that some disc $D$ contains $k$ of its zeros.
 	Then we can search for another disc $D_1$ that contains $m$ zeros of some derivative of $p(z)$, where $D_1$ depends only on disc $D$ and  parameters $n,k$ and $m$. It is also of special interest to find best possible (i.e. smallest) disc $D_1$.This problem is known as problem of Kakeya, and in general it is not solved yet. 
 	For example, in the case of the first derivative, if $k=2$ and $m=1$, then best possible option for $D_1$ is given by well-known Alexander-Kakeya theorem [4].
 	In the present paper we discuss case of $(k-1)$-th derivative of  $p(z)$, for $k > 2$, and $m=1$. In other words
 	we want to locate a single zero of $(k-1)$-th derivative of $p(z)$, provided that $k$ zeros of $p(z)$ lie in some disc $D$.

	We shall consider normalized case where $D$ is a closed unit disc. A special case when \\ $p(z) = (z^k-1)q(z)$ (i.e. when these $k$ zeros are $k$-th roots of unity) has been already studied in [5].
	
	Now, let $z_1,z_2,...,z_k$ be $n$ pairwise distinct complex numbers. We shall assume that they lie in the closed unit disc. Let $a_1,a_2,...,a_k$  be solution of the following system of $k$ linear equations:
	\begin{align*}
	a_1 z_1^m+a_2 z_2^m+...+a_k z_k ^m= \delta_{m,k-1},
	\end{align*}
	    for  $0\le m\le k-1$
	(where $\delta_{i,j}$ is the Kronecker $\delta$ symbol).

	This solution exists, since $z_i$ are pairwise distinct.
	For such $a_i$  and $z_i$ let us define $S_m$ with
\begin{align*}
 	S_m= a_1z_1^m+a_2z_2^m+...+ a_k z_k^m
 \end{align*}
	for any non-negative integer $m$.
	Our goal now is to estimate $|S_m|$.
	Let $$q(z) = (z-z_1)(z-z_2)...(z-z_k) = \sum_{i=0}^k c_i z^i .$$
	Then obviously holds
	\begin{align*}
			z_1^t q(z_1) +z_2^t q(z_2) +...+z_k^t q(z_k) = 0
	\end{align*}
	 
	for any non-negative integer $t$.
	
	Previous equality can be written in the form
	\begin{align*}
	S_{k+t}+ \quad c_{k-1}S_{k-1+t}+...+\quad c_{0}S_{t}=0
	\end{align*} 
	which give us a recurrence formula
\begin{align}
		S_{m}= -( c_{k-1}S_{m-1}+...+ c_0 S_{m-k} ).                               
\end{align}	

	Now we want to (inductively) show that
\begin{align}
	|S_{m}|\le \alpha^{m-k+1}
\end{align}
where $\alpha=\frac{1}{ \sqrt[k]{2}-1  }< \frac{k}{ \ln 2}$.  For such $\alpha$ holds $(1 + \alpha)^k - 2\alpha^k =0.$

	Obviously (2) holds for  $0\le m \le k-1$, hence we can assume $m \ge k$. Then from (1) and Viete's rules follows that
\begin{align*}
	|S_{m}|\le & |c_{k-1}||S_{m-1}|+ |c_{k-2}||S_{m-2}|+ ... +|c_0||S_{m-k}|\\
	\le& {k\choose 1 }  |S_{m-1}|+ {k\choose 2} |S_{m-2}| +...+{ k\choose k}|S_{m-k}|
	\\ \le &{k\choose 1 } \alpha^{m-k}+ {k\choose 2} \alpha^{m-k-1} +...+{ k\choose k}\alpha^{m-2k+1}
	\\ = &  \quad ((1 + \alpha)^k - \alpha^k )\alpha^{m-2k+1}.
\end{align*}
	On the other hand we have that
	$((1 + \alpha)^k - \alpha^k )\alpha^{m-2k+1} = \alpha^{m-k+1}$
	which  completes our induction.

	Now let  $a(z)=\sum_{k=0}^n a_k z^k $  and $b(z)=\sum_{k=0}^m b_k z^k$  be two complex polynomials of degree $n$ and $m$, $m\le n$.
	If $m<n$, then we shall assume that $b_{m+1} = b_{m+2} =...= b_{n}= 0$.
	
	For these polynomials we can define operator  
	$A(a,b)=\sum_{k=0}^n (-1)^k  \frac{ a_k b_{n-k}}{ {n\choose k}   }.$   
	If $n = m$ and \\ $A(a,b)=0$, then $a$ and $b$ are said to be apolar polynomials. For apolar polynomials holds classical theorem due to Grace: 
	\begin{thm}(Theorem of Grace) 
		If all zeros of a polynomials $a(z)$ are contained in some circular region  $R$, then at least one zero of $b(z)$ is contained in $R$, provided that $a(z)$ and $b(z)$ are apolar.
	\end{thm} 
 
	We recall that circular region is (open or closed) disc or half-plane or their exterior. Some generalization of the theorem of Grace can be found  in [1]-[4].
	In the case of $m\le n$ and \\ $A(a,b)=0$ this polynomials are sometimes called weakly apolar.
	The following theorem is a kind of a generalization of the theorem of Grace.
	\begin{thm} 
	  Let  $a(z)=\sum_{k=0}^n a_k z^k $  and $b(z)=\sum_{k=0}^m b_k z^k $  be two complex polynomials of degree $n$ and $m$, $m\le n$. If  $A(a,b)=0$, then polynomials $a^{(n-m)}(z)$ and $b(z)$ are apolar.
	\end{thm}
	Above theorem has been proved in [3], and follows
	from the equality
	$$A(a,b) = \frac{   (-1)^{n-m}}{n(n-1)...(n-m+1)} A( a^{(n-m)}, b ).$$
	
	Suppose now that $z_1,z_2,...,z_k$ are pairwise distinct
	zeros of a complex polynomial $p(z)$ of degree $n$. Then $A(p(z), a_i (z-z_i)^n) = 0$, for all $i$ and for any $a_i$. That implies that \\ $A (p, \sum_{i=0}^k  a_i (z-z_i )^n ) = 0$, i.e. these two polynomials are weakly apolar. Now, due to our preliminary considerations, we can choose $a_i$ such that
	\begin{align*}
		a_1 z_1^m+a_2 z_2^m+...+a_k z_k^m= \delta_{m,k-1},
	\end{align*}
	   for  $0\le m\le k-1$
	which means that  
	\begin{align*}
		t(z)=\sum_{i=1}^{k}a_i (z-z_i ) ^n= \sum_{i=k-1}^{n}{n \choose i}S_{i}(-1)^{i}z^{n-i}=\sum_{i=0}^{n-k+1}t_{i}z^{i},
	\end{align*} $t_{i}=(-1)^{n-i}{n\choose i}S_{n-i}$. 
	is a polynomial of degree $n-k+1$.  From Theorem 1 we can conclude that $p^{(k-1)} (z)$ and $t(z)$ are apolar. By theorem of Grace it
	follows that any disc that contains zeros of $t(z)$ must contain at least one zero $p^{(k-1)} (z)$. 
	So, our next task is to find disc that contains all zeros of $t(z)$. In order to do that we shall use the following well-known citerion:
	\begin{thm}
  Let  $a(z)=\sum_{k=0}^n a_k z^k  $ be a complex polynomials of degree $n$. Then all its zeros lie in the disc $|z|\le M$, where
	$$M =  2  \max_{i} \{ \sqrt[n-i]{ |\frac{a_i  }{a_n}  |  }, i=0,1,...,n-1       \}.$$ 	\end{thm}

	In the case of $t(z)$ it is easy to see that 
	\begin{align*}
	M = 2| \frac{ t_{n-k}}{t_{n-k+1}} |= 2 \frac{n-k+1}{k} |S_{k}|\le  2\frac{n-k+1}{ k} \alpha
	= 2 \frac{n-k+1}{ k(\sqrt[k]{2}  -1)    }   < \frac{2(n-k+1)} {ln2}.
	\end{align*}
 
	Therefore we have proved:
	\begin{thm}
	Let $p(z)$ be a complex polynomial of degree $n$,
	having $k$ its pairwise distinct zeros in the closed unit disc. Then at least one zero of $p^{(k-1)}(z)$ lies in the disc $|z|<\frac{2(n-k+1)} {ln2}$.
	\end{thm} 
	In the case when these zeros are not pairwise distinct 
	we can also obtain above theorem for a disc $|z|\le  \frac{2(n-k+1)} {ln2}$.
	Naimely, if we suppose that these $k$ zeros are not necessarily 
	distinct and the above theorem is not true, that means that
	$p^{(k-1)}(z)$ has no zero in the disc $|z|\le \frac{2(n-k+1)} {ln2}$. Then we can "very slightly" shift these zeros inside unit disc, making them
	pairwise distinct and keeping also zeros of  $p^{(k-1)}(z)$ outside disc
	$|z|\le \frac{2(n-k+1)} {ln2}$. But this situation is in contradiction with Theorem 3. Therefore we have proved:

\begin{thm}
Let $p(z)$ be a complex polynomial of degree $n$,
having $k$ its  zeros in the closed unit disc. Then at least one zero of $p^{(k-1)}(z)$ lies in the disc $|z|\le \frac{2(n-k+1)} {ln2}$.
\end{thm} 
	If $k=n$ ( i.e. if we want to locate a zero of $p^{(n-1)}(z)$ provided that all zeros of $p(z)$ lie in the closed unit disc) then the only zero of $p^{(n-1)}(z)$ lies in the closed unit disc for any polynomial $p(z)$ having all its zeros in it. Since disc $|z|\le \frac{2(n-k+1)} {ln2}$   contains the closed unit disc that means that our estimation is of no use for $k =n$,  but it has a sense for other values of  $k$.
	
	At the end let us not that with above approach we can locate a zero of not only $p^{(n-1)}(z)$, but also of $p^{(i)}(z)$, for
	$1\le i<k-1$. In that case starting point would be to replace $\delta_{m,k-1}$ with $\delta_{m,i}$ (in the beginning part), but we find case $i=k-1$ as the most interesting.

\end{document}